\documentclass{ita}

\usepackage{amssymb}
\usepackage{url}
\usepackage{gastex}

\theoremstyle{definition}
\newtheorem{fact}[thrm]{Fact}

\def\calA{\mathcal{A}}
\def\calB{\mathcal{B}}
\def\calC{\mathcal{C}}

\def\form#1{\mathsf{#1}}

\def\card{\form{card}}

\def\rank{\form{rank}}

\def\diam{\form{diam}}

\def\inv{^{-1}}
\let\phi\varphi

\def\leff{\le_{\sf ff}}

\def\mapright#1{\mathop{\longrightarrow}\limits^{#1}}
\def\mapleft#1{\mathop{\longleftarrow}\limits^{#1}}

\def\expstep{\longrightarrow_{{\sf exp}}}
\def\estep{\longrightarrow_{{\sf e}}}
\def\istep{\longrightarrow_{{\sf i}}}
\def\restep{\longrightarrow_{{\sf re}}}

\begin{document}

\title{On an algorithm to decide whether a free group is a free factor
of another}
\thanks{The first author acknowledges support from C.M.U.P., financed
by F.C.T. (Portugal) through the programmes POCTI and POSI, with
national and European Community structural funds.  Both authors
acknowledge support from the European Science Foundation program
AutoMathA.}

\author{Pedro V. Silva}
\address{Centro de Matem\'{a}tica, Faculdade de Ci\^{e}ncias --
Universidade do Porto -- R. Campo Alegre 687 -- 4169-007 Porto,
Portugal; email: \url{pvsilva@fc.up.pt}}
\author{Pascal Weil}
\address{LaBRI, CNRS -- 351 cours de la Lib\'eration -- 33405 Talence
Cedex -- France; email: \url{pascal.weil@labri.fr}}
\subjclass{20E05,05C25}
\keywords{combinatorial group theory, free groups, free factors, 
inverse automata, algorithms}

\begin{abstract}
    We revisit the problem of deciding whether a finitely generated
    subgroup $H$ is a free factor of a given free group $F$.  Known
    algorithms solve this problem in time polynomial in the sum of the
    lengths of the generators of $H$ and exponential in the rank of
    $F$.  We show that the latter dependency can be made exponential
    in the rank difference $\rank(F) - \rank(H)$, which often makes a
    significant change.
\end{abstract}

\begin{resume}
    Nous revenons sur la question de d\'ecider si un sous-groupe
    finiment engendr\'e $H$ est facteur libre d'un groupe libre
    donn\'e $F$.  On trouve dans la litt\'erature des algorithmes qui
    r\'esolvent ce probl\`eme en temps polynomial en la somme des
    longueurs des g\'en\'erateurs de $H$, et exponentiel en le rang de
    $F$.  Nous montrons que l'on peut remplacer la d\'ependance
    exponentielle en $\rank(F)$ par une d\'ependance exponentielle en
    la diff\'erence $\rank(F) - \rank(H)$, ce qui change souvent les
    choses de fa\c{c}on consid\'erable.
\end{resume}

\maketitle

%%%%%%%%%%%%%%%%%%%%%%%%%%%%%%%

The combinatorial aspects of group theory have attracted the attention
of theoretical computer scientists for a long time, and for a variety
of reasons.  There is no need to recall the importance of the concept
of monoid (e.g. free, finite) in the theory of automata since the
foundational results of Sch\"utzenberger and Eilenberg (see the books
\cite{PinBook,Almeida}), and groups form a special case of monoids
that sometimes play an important role in purely language- and
monoid-theoretic questions (e.g. the type II conjecture, see
\cite{HMPR} for a survey).  Algorithmic questions (the word problem,
the conjugacy problem, \etc.)  have been very influential in group
theory throughout the 20th century, starting from the work of Dehn,
and specialists of combinatorics on words find a particular interest
in the challenges posed by the analogous combinatorics of the free
group.  Recent work on important non-commutative groups like the
Thompson group and the so-called automata groups strongly relies on
the formalism of finite state automata (see
\cite{cleary,KambitesSS,SilvaSteinberg} for recent examples).  It is
already a classical result that these automata-theoretic and
combinatorial points of view converge in the (admittedly simpler)
study of the subgroups of free groups, this is central in the
algorithmic problem tackled here and is discussed in detail in the
first part of this paper.

Let us also mention another reason for the recent multiplication of
research projects on the boundary between computer science and
combinatorial group theory.  Public-key cryptography relies heavily on
group theory: sometimes finite groups such as the groups of units in
modular arithmetic, or the groups of rational points on elliptic
curves over finite fields, sometimes infinite non-commutative groups
like the braid groups (see for instance
\cite{sidelnikov,anshel,dehornoy} and many others).  At any rate, the
design of more robust cryptographic schemes and the attack of such
schemes rely on a deeper understanding of the combinatorial and
algorithmic properties of non-commutative groups.

As mentioned above, the combinatorial and algorithmic problems
concerning free groups are of special interest.  Free groups are
archetypal groups, whose structure is far from being totally
elucidated, and the efficient solution of standard problems in their
context can shed some light on the possible solution of the same
problems in more complex groups.  Moreover, the solution of
algorithmic problems in free groups may be more attainable since we
can use the resources of combinatorics on words and automata theory.

We now present the specific algorithmic problem addressed in this
paper.  For the classical facts about free groups recorded below
without a reference, we refer the reader to the book by Lyndon and
Schupp \cite{LS}.  It is well-known that the minimal sets of
generators, or \textit{bases}, of a free group $F$ all have the same
cardinality, called the \textit{rank} of $F$.  Moreover, if $F$ has
finite rank $r$, every $r$-element generating set of $F$ is a basis,
see \cite[Prop.  I.3.5]{LS}.  In this paper, we consider only finite
rank free groups.

Let $H$ be a subgroup of a free group $F$, written $H\le F$. Then $H$
itself is a free group whose rank may be greater than the rank of $F$.
We say that $H$ is a \textit{free factor} of $F$, written $H\leff F$,
if there exist bases $B$ of $H$ and $A$ of $F$ such that $B\subseteq
A$ (free factors can be defined in all groups by a universal property,
but the operational definition given here is sufficient for the
purpose of this study).  It is well known that one can decide whether
a given finite rank subgroup $H\le F$ is a free factor of $F$, but the
known algorithms have a rather high time complexity.  More precisely,
the best of these algorithms require time that is polynomial in the
size of $H$ and exponential in the rank of $F$, see
Section~\ref{whitehead} below for the details.  Here, the size of $H$
is taken to be the sum of the lengths of a finite set of generators of
$H$ in $F$.

We propose a new algorithm to decide whether a given finitely
generated subgroup $H$ is a free factor of the free group $F$, which
is polynomial in the size of $H$ and exponential in the rank
difference between $F$ and $H$.  In many instances, this represents a
substantial advantage over exponential dependency in the rank of $F$.

Our algorithm relies essentially on a careful analysis of the
construction of the \textit{graph representation} of $H$.  More
precisely, once a basis $A$ of the ambient free group $F$ is fixed,
there is a natural and elegant representation of the finitely
generated subgroups of $F$ by $A$-labeled graphs (or inverse
automata).  This construction --- a graphical representation of ideas
that go back to the early part of the twentieth century \cite[Chap.
11]{Rotman} --- was made explicit by Serre \cite{Serre} and Stallings
\cite{Stallings}.  It has been used to great profit by many authors
since the late 1970s, see
\cite{MargolisMeakin,Ventura,MargolisSapirWeil,KapovichMiasnikov} for
recent examples.  Given a finite set of generators of $H$ (as reduced
words over the alphabet $A\cup A\inv$), the graph representation of
$H$ can be effectively constructed (see \cite{Stallings},
\cite{MargolisSapirWeil}, \etc).  The number of vertices and edges of
this graph is bounded above by $\ell$, the sum of the lengths of a set
of generators of $H$, and the whole representation can be computed in
time at most $O(\ell^2)$ (in fact, in time $O(\ell\log^*\ell)$
according to a recent announcement\footnote{For a positive integer
$n$, $\log^{*}(n)$ is the least integer $k$ such that the $k$-th
iterate of the $\log$ function of $n$ is at most 1.  The growth of
$\log^{*}(n)$ is so slow that it can be considered a constant for all
practical purposes\dots} by Touikan~\cite{Touikan}).  We discuss this
representation in more detail in Section~\ref{representation} below,
and we show in Sections~\ref{careful} and \ref{sec deciding} how to
use it to decide more efficiently the free factor relation.

It is interesting to note that our algorithm is the first to be
expressed entirely in terms of the graph representation of $H$.  Let
us also emphasize that we do not claim that our algorithm is optimal.
It is an open question whether one can decide the free factor relation
$H \leff K$ in time polynomial both in the size of $H$ and in the
size, or the rank of $K$.

%%%%%%%%%%%%%%%%%%%%%%%%%%%%%%%
\section{Background}

If $A$ is a basis of a free group $F$, we often write $F = F(A)$ and
we represent the elements of $F$ as reduced words over the alphabet
$A$.  More precisely, we consider the set of all words on the
symmetrized alphabet $\tilde A = A \cup A\inv$, where $A\inv = \{a\inv
\mid a\in A\}$ is a set that is disjoint from $A$, equipped with an
explicit bijection with $A$, namely $a\mapsto a\inv$.  It is customary
to extend the mapping $u \mapsto u\inv$ to all words $u\in \tilde
A^{*}$ by letting $(a\inv)\inv = a$ for each $a\in A$, $1\inv = 1$
(where 1 denotes the empty word) and $(a_{1}a_{2}\cdots a_{n})\inv =
a_{n}\inv \cdots a_{2}\inv a_{1}\inv$ for all $a_{1},\ldots,a_{n}\in
\tilde A$.  A word in $\tilde A^{*}$ is \textit{reduced} if it
contains no factor of the form $aa\inv$ or $a\inv a$ with $a\in A$,
and it is well known that $F$ can be identified with the set of
reduced words over $A$.  We denote by $\rho$ the map that assigns to
each word $u$ the corresponding reduced word $u\rho \in F(A)$,
obtained by iteratively deleting all factors of the form $aa\inv$ or
$a\inv a$ ($a\in A$).

%%%%%%%%%%%%%%%%%%%%%%%%%%%%%%%
\subsection{On inverse automata}

We describe the main tool for the representation of subgroups of free
groups in terms of automata (see \cite{PerrinHB}). Readers less
familiar with this terminology may think of automata as edge-labeled
directed graphs.

An \textit{automaton} on alphabet $A$ is a triple of the form $\calA =
(Q, q_{0}, E)$ where $Q$ is a finite set called the \textit{state}
set, $q_{0}\in Q$ is the \textit{initial state}, and $E \subseteq Q
\times A \times Q$ is the set of \textit{edges}, or
\textit{transitions}.  A transition $(p,a,q)$ is said to be from state
$p$, to state $q$, with label $a$.  The \textit{label of a path} in
$\calA$ (a finite sequence of consecutive transitions) is the sequence
of the labels of its transitions, a word on alphabet $A$, that is, an
element of the free monoid $A^*$.  We write $p \mapright u q$ if there
is a path from state $p$ to state $q$ with label $u$.  The language
accepted by $\calA$ is the set $L(\calA)$ of all words in $A^*$ which
label a path in $\calA$ from $q_{0}$ to $q_{0}$.

This definition of automata leads naturally to the definition of a
\textit{homomorphism} $\phi$ from an automaton $\calA = (Q,q_{0},E)$
to an automaton $\calA' = (Q',q'_{0},E')$ (over the same alphabet
$A$): $\phi$ is a mapping from $Q$ to $Q'$ such that $\phi(q_{0}) =
q'_{0}$, and such that whenever $(p,a,q)\in E$, we also have
$(\phi(p),a,\phi(q)) \in E'$.  The homomorphism $\phi$ is an
\textit{isomorphism} if it is a bijection and if $\phi\inv$ is also a
homomorphism.

The automaton $\calA$ is called \textit{deterministic} if no two
distinct edges with the same initial state bear the same label, that
is,
$$(p,a,q), (p,a,q') \in E \Longrightarrow q = q'.$$
The automaton is called \textit{trim} if every state $q\in Q$ lies in
some path from $q_{0}$ to $q_{0}$.

In the sequel, we consider automata where the alphabet is symmetrized,
that is, the alphabet is of the form $\tilde A = A \cup A\inv$.  We
say that $\calA$ is \textit{dual} if for each $a\in A$, there is an
$a$-labeled edge from state $p$ to state $q$ if and only if there is
an $a\inv$-labeled edge from $q$ to $p$,
$$(p,a,q) \in E \Longleftrightarrow (q,a\inv,p)\in E.$$
Let us immediately record the following fact.

\begin{fact}\label{fact 11}
    Let $\calA$ be a deterministic dual automaton.  If a word $u$
    labels a path in $\calA$ from state $p$ to state $q$, then so does
    the corresponding reduced word $u\rho$.
    Moreover $L(\calA)$ is a submonoid of $\tilde A^*$ and
    $L(\calA)\rho$ is a subgroup of $F(A)$.
\end{fact}    

Now let $\calA = (Q,q_{0},E)$ be a trim dual automaton and let $p,q\in
Q$ be states of $\calA$. If $w = a_{1}\cdots a_{n}\in \tilde A^*$
is a non-empty word, the expansion of $\calA$ by $(p,w,q)$ is the
automaton obtained from $\calA$ by adding $n - 1$ vertices
$q_{1},\ldots,q_{n-1}$ and $2n$ edges
$$p \mapright{a_1} q_1 \mapright{a_2} \ldots \mapright{a_{n-1}}
q_{n-1} \mapright{a_{n}} q$$
and
$$q \mapright{a_{n}\inv} q_{n-1} \mapright{a_{n-1}\inv} \ldots
\mapright{a_2\inv} 
q_{1} \mapright{a_1\inv} p.$$
Note that this automaton is still trim and dual. Moreover, if $p = q 
= q_{0}$, then we observe the following.

\begin{prpstn}\label{H(A)}
    Let $\calA = (Q,q_{0},E)$ be a trim dual automaton, let $H =
    L(\calA)\rho$ and let $w$ be a non-empty word. If $\calB$ is the
    expansion of $\calA$ by $(q_{0},w,q_{0})$, then $L(\calB)\rho$ is
    the subgroup generated by $H$ and $w\rho$, that is, $L(\calB)\rho =
    \langle H, w\rangle$.
\end{prpstn}

\begin{proof}
Let $\calC$ be the dual automaton consisting of the state $q_{0}$ and
the states and edges added to $\calA$ in the expansion. It is
immediate that $L(\calC)\rho$ is the subgroup of $F(A)$ generated by 
$w\rho$.

If $u\in L(\calB)$, we can factor a path $q_{0} \mapright u q_{0}$
according to the successive visits of state $q_{0}$. The resulting
factorization of $u$ makes it clear that $u$ is a product of elements
of $L(\calA)$ and $L(\calC)$. Thus, $L(\calB)$ is the submonoid
generated by $L(\calA) \cup L(\calC)$, and $L(\calB)\rho$ is the
subgroup generated by $L(\calA)\rho$ and $w\rho$. This concludes the
proof.
\end{proof}

%%%%%%%%%%%%%%%%%%%%%%%%%%%%%%%
\subsection{Reduced inverse automata}\label{representation}

The automaton $\calA$ is called \textit{inverse} if it is
deterministic, trim and dual.  It is \textit{reduced} if every state
$q\in Q$ lies in some path from $q_{0}$ to $q_{0}$, labeled by a
(possibly empty) reduced word.  We note the following result, a cousin
of \cite[Thm 1.16]{Stephen}.

\begin{prpstn}\label{prop unique}
    If $\calA$ and $\calB$ are reduced inverse automata such that
    $L(\calA)\rho = L(\calB)\rho$, then $\calA$ and $\calB$ are
    isomorphic.
\end{prpstn}

\begin{proof}
Let $\calA = (Q,q_{0},E)$ and $\calB = (P,p_{0},D)$ be reduced
inverse automata such that $L(\calA)\rho = L(\calB)\rho$. We
construct an isomorphism $\phi$ between $\calA$ and $\calB$ as follows.
We first let $\phi(q_{0}) = p_{0}$.

Let $q\in Q$.  Since $\calA$ is reduced, there exist reduced words $u$
and $v$ such that the word $uv$ is reduced, $q_{0} \mapright{u} q$ and
$q\mapright{v} q_{0}$.  Then $uv\in L(\calA)\rho$, so $uv \in
L(\calB)\rho$, and hence $uv\in L(\calB)$ by Fact~\ref{fact 11}.  Thus
$uv$ labels a path in $\calB$ from $p_{0}$ to $p_{0}$, and we let
$\phi(q)$ be the unique state in $P$ such that $p_{0} \mapright{u}
\phi(q) \mapright{v} p_{0}$.

We first verify that $\phi$ is well defined.  Suppose that $uv$ and
$u'v'$ are reduced words such that $q_{0} \mapright{u} q \mapright{v}
q_{0}$ and $q_{0} \mapright{u'} q \mapright{v'} q_{0}$ in $\calA$.  We
want to show that if $p_{0} \mapright{u} p \mapright v p_{0}$ and
$p_{0} \mapright{u'} p' \mapright{v'} p_{0}$ in $\calB$, then $p =
p'$.  We note that $u'v$ labels a path from $q_{0}$ to $q_{0}$ in
$\calA$.  If $u'v$ is a reduced word, then by the same reasoning as
above, $u'v$ labels a path in $\calB$ from $p_{0}$ to $p_{0}$, say,
$p_{0} \mapright{u'} p'' \mapright v p_{0}$ and the deterministic
property of $\calB$ implies that $p' = p'' = p$.

If $u'v$ is not reduced, and $a$ is the first letter of $v$, then the
last letter of $u'$ is $a\inv$ while the last letter of $u$ is not
$a\inv$. Therefore $u'u\inv$ is reduced, $u'u\inv \in L(\calA)$ and
again, there is a path in $\calB$ of the form $p_{0} \mapright{u'} p''
\mapright{u\inv} p_{0}$. By determinism, it follows that $p' = p'' =
p$.

This shows that $\phi$ is well defined.  A dual construction yields a
well-defined mapping $\psi$ from $P$ to $Q$ such that, whenever $p_{0}
\mapright u p \mapright v p_{0}$ in $\calB$ and $uv$ is a reduced
word, then $q_{0} \mapright u \psi(p) \mapright v q_{0}$ in $\calA$.
Using the determinism of $\calA$ and $\calB$, it is now immediate that
$\psi\circ \phi$ is the identity on $Q$ and $\phi\circ \psi$ is the
identity on $P$.

There remains to verify that $\phi$ and $\phi\inv$ are homomorphisms.
The case of $\phi\inv$ is dual of that of $\phi$ and we treat only the
latter.  That is, we want to show that if $(q,a,q')$ is a transition
in $\calA$, then $(\phi(q),a,\phi(q'))$ is a transition in $\calB$.
Let $uv$ and $u'v'$ be reduced words such that $q_{0} \mapright u q
\mapright v q_{0}$ and $q_{0} \mapright{u'} q' \mapright{v'} q_{0}$.
In particular, we have $p_{0} \mapright u
\phi(q) \mapright v p_{0}$ and $p_{0} \mapright{u'} \phi(q')
\mapright{v'} p_{0}$ in $\calB$.

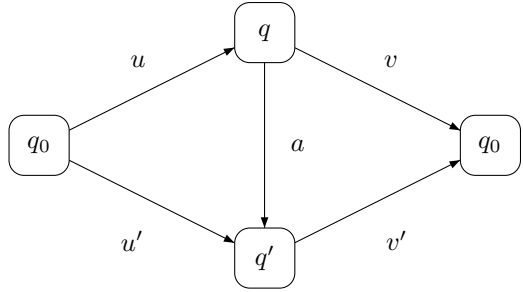
\begin{center}
    \begin{picture}(60,38)(0,-34)

	\node[NLangle=0.0,Nmr=2.0](n0)(0.0,-15.0){$q_0$}

\node[NLangle=0.0,Nmr=2.0](n3)(30.0,0.0){$q$}

\node[NLangle=0.0,Nmr=2.0](n4)(30.0,-30.0){$q'$}

\node[NLangle=0.0,Nmr=2.0](n5)(60.0,-15.0){$q_0$}

\drawedge[ELdist=3.0](n0,n3){$u$}

\drawedge[ELdist=3.0](n3,n5){$v$}

\drawedge[ELdist=3.5](n3,n4){$a$}

\drawedge[ELside=r,ELdist=4.0](n0,n4){$u'$}

\drawedge[ELside=r,ELdist=4.0](n4,n5){$v'$}

\end{picture}
\end{center}

If $uav'$ is reduced, then in $\calB$, there is a path from $p_{0}$ to
$p_{0}$ labeled $uav'$, and by determinism, there is a transition
$(\phi(q),a,\phi(q'))$. If $uav'$ is not reduced, then either $ua$ is
not reduced or $av'$ is not reduced. If $ua$ is not reduced, then $u =
u_{1}a\inv$ and by determinism, $q_{0} \mapright{u_{1}} q'$. As in the
first part of the proof, it follows that at least one of $u_{1}v'$ and
$u_{1}{u'}\inv$ is reduced, so $p_{0} \mapright{u_{1}} \phi(q')$ in
$\calB$ and hence there is a transition $(\phi(q),a,\phi(q'))$.
The case where $av'$ is not reduced is handled symmetrically, and
this concludes the proof.
\end{proof}

Let $H$ be a subgroup of $F(A)$.  Say that an automaton $\calA$
\textit{on alphabet $A$ represents $H$} if $\calA$ is reduced and
inverse and if $L(\calA)\rho = H$.  Proposition~\ref{prop unique}
shows that there exists at most one such automaton, and we denote it
by $\Gamma_{A}(H)$ if it exists.  We now discuss the existence and the
construction of $\Gamma_{A}(H)$ when $H$ is finitely generated.  (As
it turns out, $\Gamma_{A}(H)$ always exists, but our interest in this
paper is restricted to the finite rank case.)

Let $\calA$ be an automaton and let $p, q$ be distinct states of
$\calA$. The automaton obtained from $\calA$ by identifying states $p$
and $q$ is constructed as follows: its state set is
$Q\setminus\{p,q\}\cup \{n\}$, where $n$ is a new state; its initial
state is $q_{0}$ (or $n$ if $p$ or $q$ is equal to $q_{0}$); and its
set of transitions is obtained from $E$ by replacing everywhere $p$ and
$q$ by $n$. If $\calA$ is trim or dual, then so is the automaton
obtained from $\calA$ by identifying a pair of states.

Now let $\calA$ be a dual automaton.  If $\calA$ is not deterministic,
there exist transitions $(r,a,p)$ and $(r,a,q)$ with $p\ne q$ and
$a\in \tilde A$.  Identifying $p$ and $q$ yields a new dual automaton
$\calB$, and we say that $\calB$ is obtained from $\calA$ by an
\textit{elementary reduction of type 1}.

\begin{fact}\label{reduction1}
    Let $\calA$ be a dual automaton and let $\calB$ be obtained from
    $\calA$ by an elementary reduction of type 1. Then $L(\calA)\rho =
    L(\calB)\rho$.
\end{fact}

\begin{proof}
It is easily seen that $L(\calA) \subseteq L(\calB)$. For the converse,
we use the notation given above: in $\calB$, the states $p$ and $q$ of
$\calA$ are replaced with a new state $n$. Let $u\in L(\calB)$. Then
there exists a path labeled $u$ from the initial state of $\calB$ (say,
$q_{0}$) to itself. If that path does not visit state $n$, then $u$
also labels a path from $q_{0}$ to itself in $\calA$ and hence $u\in
L(\calA)$.

If that path does visit state $n$, we consider the
factorization of $u$ given by the passage of that path through $n$: we
have $u = u_{0}u_{1}\cdots u_{r}$, $r\ge 1$ and
$$q_{0} \mapright{u_{0}} n \mapright{u_{1}} n \cdots n \mapright{u_{r}}
q_{0}.$$
It follows that in $\calA$, $u_{i}$-labelled paths exist, with end
states $p$ or $q$ (or $q_{0}$).  Then one of $u_{0}$ and $u_{0}a\inv
a$ labels a path in $\calA$ from $q_{0}$ to $q$.  Similarly, one of
$u_{r}$ and $a\inv au_{r}$ labels a path from $q$ to $q_{0}$ (making
due allowance if $p$ or $q$ is equal to $q_{0}$).  And for each $1\le
i \le r$, one of $u_{i}$, $a\inv au_{i}$, $u_{i}a\inv a$ and $a\inv
au_{i}a\inv a$ labels a path in $\calA$ from $q$ to $q$.  Therefore,
there exists a path in $\calA$ of the form $q_{0} \mapright v q_{0}$
such that $u\rho = v\rho$, which concludes the proof.
\end{proof}

Now assume that $\calA$ is a deterministic dual automaton. 

\begin{fact}\label{reduction11}
    Let $\calA$ be an inverse automaton.  Then $\calA$ is non-reduced
    if and only if there exist states $q\ne q_{0}$ and $p$, and a
    letter $a\in \tilde A$ such that the only transitions of $\calA$
    involving $q$ are $(p,a,q)$ and $(q,a\inv,p)$.
    
    In graph-theoretic terms, this means that $\calA$ is reduced if
    and only if no vertex of $\calA$ has degree one (more precisely: 
    no vertex is adjacent to a single $A$-labeled edge), except possibly
    $q_{0}$.
\end{fact}

\begin{proof}
By definition, $\calA$ is not reduced if and only if there exists a
state $q$ that does not lie on any path from $q_{0}$ to itself,
labeled by a reduced word.  We first observe that the state $q$ cannot
be equal to $q_{0}$ since the empty word is reduced, and labels a path
from $q_{0}$ to itself.  Suppose now that there exist reduced words
$u,v$ with distinct last letters, labeling paths from $q_{0}$ to $q$:
then $uv\inv$ is a reduced word, labeling a path from $q_{0}$ to
itself and visiting $q$.  On the other hand, if every reduced word
labeling a path from $q_{0}$ to $q$ ends with, say, letter $a\in
\tilde A$, then every path from $q_{0}$ to itself visiting $q$ has a
non-reduced label.  Thus $\calA$ is not reduced if and only if there
exists a state $q\ne q_{0}$ and every reduced word labeling a path
from $q_{0}$ to $q$ ends with the same letter.  By determinism, this
is equivalent to the existence of another state $p$ such that the
transitions involving $q$ are $(p,a,q)$ and $(q,a\inv,p)$.
\end{proof}

Let $\calA$ be inverse and not reduced, and let $a,p,q$ be as in
Fact~\ref{reduction11}.  If $\calB$ is obtained from $\calA$ by
omitting state $q$ and the transitions involving it, we observe that
$\calB$ is again an inverse automaton, and we say that $\calB$ is
obtained from $\calA$ by an \textit{elementary reduction of type 2}.

\begin{fact}\label{reduction2}
    Let $\calA$ be an inverse automaton and let $\calB$ be obtained
    from $\calA$ by an elementary reduction of type 2.  Then
    $L(\calA)\rho = L(\calB)\rho$.
\end{fact}

\begin{proof}
Let $a$ be a letter and let $p,q$ be states of $\calA$ as in
Fact~\ref{reduction11}, and let us assume that $\calB$ is obtained
from $\calA$ by omitting state $q$ and the transitions involving it.
It is easily seen that $L(\calB) \subseteq L(\calA)$.  Conversely, let
$u\in L(\calA)$.  By Fact~\ref{fact 11}, $u\rho\in L(\calA)$.  Now
Fact~\ref{reduction11} shows that the path $q_{0} \mapright u\rho
q_{0}$ in $\calA$ cannot visit state $q$, since $u\rho$ is a reduced
word.  It follows that this path is also a path in $\calB$, $u\rho\in
L(\calB)$ and hence $u\rho\in L(\calB)\rho$.
\end{proof}

Let $\calA$ be a trim, dual automaton, and let $\calB$ be an automaton
obtained by iteratively performing elementary reductions, first of
type 1 until the automaton is inverse, and then of type 2 until none
is possible.  Then $\calB$ is a reduced inverse automaton, we write
$\calB = \calA\rho$ and we say that $\calB$ is obtained from $\calA$
by reduction.  Moreover, Facts~\ref{reduction1} and \ref{reduction2}
show that $L(\calA)\rho = L(\calB)\rho$.

This leads directly to the well-known algorithm to construct a reduced
inverse automaton representing a given finitely generated subgroup
$H$.  Let $h_{1},\ldots,h_{n}$ be generators of $H$, and let us
consider the automaton obtained from the trivial automaton (one vertex
$q_{0}$, no transitions) by performing successively expansions by
$(q_{0}, h_{i}, q_{0})$ ($1\le i\le n$) and then reducing the
automaton.  It follows from Propositions~\ref{H(A)} and~\ref{prop
unique} that the resulting automaton is $\Gamma_{A}(H)$.  Note that it
does not matter which set of generators of $H$ was used, nor in which
order the elementary reductions were performed.

\begin{rmrk}
    This construction of $\Gamma_{A}(H)$ is well known, and can be
    described in many different ways, notably in terms of immersions
    over the bouquet of circles (Stallings \cite{Stallings}) or of
    closed inverse submonoids of a free inverse monoid (Margolis and
    Meakin \cite{MargolisMeakin}).
\end{rmrk}    

\begin{fact}\label{end of sec 12}
    There is a well-known converse to the above construction: if
    $\calA$ is a reduced inverse automaton and $H = L(\calA)\rho$,
    then $H$ has finite rank and a basis for $H$ can be computed as
    follows (see Stallings \cite{Stallings}).  Given a spanning tree
    $T$ of the (graph underlying the) automaton $\calA$, for each
    state $p$, let $u_{p}$ be the reduced word labeling a path from
    $q_{0}$ to $p$ inside the tree $T$.  For each transition $e =
    (p,a,q)$, let $b_{e} = u_{p}au_{q}\inv$: then a basis of $H$
    consists of the elements $b_{e}$, where $e$ runs over the
    transitions $e = (p,a,q)$ not in $T$ and such that $a\in A$.
    
    We note that, given a finite set $h_{1},\ldots,h_{n}$ of elements
    of $F(A)$ with total length $\ell = \sum_{i}|h_{i}|$, one can
    construct $\Gamma_{A}(H)$ in time at most $O(\ell^2)$ and
    $\Gamma_{A}(H)$ has $v \le \ell - n + 1$ states.  Moreover,
    finding a basis of $H$ can be done in time at most $O(v^2)$ ($O(v
    \log^*v)$ according to Touikan's announcement~\cite{Touikan}), and
    the rank of $H$ is equal to $e-v+1$, where $e$ is the number of
    edges in $\Gamma_{A}(H)$.
\end{fact}

%%%%%%%%%%%%%%%%%%%%%%%%%%%%%%%
\subsection{On the complexity of Whitehead and other
algorithms}\label{whitehead}

It is well known that one can decide, given a subgroup $H$ of a finite
rank free group $F$, whether $H$ is a free factor of $F$. We briefly
describe here the main known algorithms and discuss their complexity.

Let $H$ be a finitely generated subgroup of a free group $F$ of rank
$r$, with basis $A$. Let $h_{1},\ldots,h_{n}$ be a generating set of
$H$. By the results summarized in Fact~\ref{end of sec 12}, up to a
quadratic time computation, we may assume that $h_{1},\ldots,h_{n}$ is
a basis of $H$. Let $\ell = |h_{1}| + \cdots + |h_{n}|$ be the total
length of the tuple $(h_{i})_{i}$, and let $d = r - n$ be the rank
difference between $F$ and $H$ -- which we assume to be positive, since
$H$ can be a proper free factor of $F$ only if $n < r$.

Federer and J\'onsson (see \cite[Prop. I.2.26]{LS}) gave the following
observation and decision procedure: $H$ is a free factor of $F$ if and
only if there exist $d$ words $h_{n+1},\ldots,h_{r}$, each of length at
most $\max\{|h_{i}| \mid 1\le i\le n\}$, such that $h_{1},\ldots,h_{r}$
generate the whole of $F$. The resulting algorithm requires testing
every suitable $d$-tuple of reduced words on alphabet $A$. Each of
these tests (does a certain $r$-tuple of words generate $F$?) takes
time polynomial in the total length of the $r$-tuple, and hence in
$d\ell$. However, the number of tests is $O(r^{d\ell})$, which is
exponential in $\ell$ and $d$.

This approach leads to the following.

\begin{fact}
    Deciding whether $H\leff K$ is in $NP$, with respect to $d\ell$.
\end{fact}

\begin{proof}
To verify that $H\leff K$, we need to guess $d$ words of length at most
$\ell$, and verify that together with $H$, they generate $F$, which can
be done in $O((d\ell)^2)$.
\end{proof}

Another approach is based on the use of Whitehead automorphisms.  We
refer the readers to \cite[Sec.  I.4]{LS} for the definition of these
automorphisms, it suffices to note here that the set $W$ of Whitehead
automorphisms of $F$ which do not preserve length, has exponential
cardinality (in terms of $r$).  A result of Whitehead \cite[Prop.
I.4.24]{LS} shows the following: if there exists an automorphism
$\phi$ of $F$ such that the total length of $(\phi(h_{i}))_{i}$ is
strictly less than $\ell$, then there exists such an automorphism in
$W$.  In particular, an algorithm to compute the minimum total length
of an automorphic image of the tuple $(h_{i})_{i}$ consists in
repeatedly applying the following step: try every automorphism
$\psi\in W$ until the total length of $(\psi(h_{i}))_{i}$ is strictly
less than the total length of $(h_{i})_{i}$; if such a $\psi$ exists,
replace $(h_{i})_{i}$ by $(\psi(h_{i}))_{i}$; otherwise, stop and
output the total length of $(h_{i})_{i}$.

This applies to the decision of the free factor relation since $H\leff
F$ if and only if there exists an automorphism $\phi$ mapping a basis
of $H$ to a subset of $A$.  Thus an algorithm consists in first
computing a basis of $H$, and assuming that $h_{1},\ldots,h_{n}$ is a
basis, verifying whether the minimum total length of
$(\phi(h_{i}))_{i}$ when $\phi$ runs over the automorphisms of $F$, is
exactly $n$.  This algorithm may require $O((\ell-n)\card(W))$ steps,
each of which consists in computing the image of a tuple of length at
most $\ell$ under an automorphism, and hence has complexity $O(\ell)$.
Thus the time complexity of this algorithm is $O(\ell^2\ \card(W))$,
which is quadratic in $\ell$ and exponential in $r$.

A variant of this algorithm was established by Gersten \cite{Gersten},
who showed that a similar method applies to find the minimum size
(number of vertices) of $\Gamma_{A}(\phi(H))$, when $\phi$ runs over
the automorphisms of $F(A)$. It is clear that $H$ is a free factor of
$F(A)$ if and only if there exists an automorphism $\phi$ such that
$\Gamma_{A}(\phi(H))$ has a single vertex. The time complexity is
computed as above, where the number of vertices of $\Gamma_{A}(H)$ is
substituted for the total length of a basis for $H$. As noted earlier,
this number of vertices is usually substantially smaller than the total
length of a basis, but the two values are linearly dependent, so the
order of magnitude of the time complexity is not modified, notably the
exponential dependence in $r$.

\begin{rmrk}
    The discussion of Whitehead's algorithm above concerns only the
    so-called \textit{easy part} of the algorithm (see for instance
    Kapovich, Miasnikov and Shpilrain \cite{KSS}).  Results by
    Miasnikov and Shpilrain \cite{MS}, Khan \cite{Kh}, and most
    recently by Donghi Lee \cite{DL} on the possible polynomial
    complexity of the \textit{hard part} of the algorithm also
    consider the rank of the ambient free group as a constant, and do
    not question the actual exponential dependence in that parameter.
\end{rmrk}    

%%%%%%%%%%%%%%%%%%%%%%%%%%%%%%%
\section{A careful look at the expansions and reductions of inverse
automata}\label{careful}

Let $\calA$ be a reduced inverse automaton.

Let $\calB$ be obtained from $\calA$ by performing an expansion, say by
$(p,w,q)$, and then reducing the resulting automaton. In this
situation, we write $\calA \expstep^{(p,w,q)} \calB$, or simply $\calA
\expstep \calB$. We distinguish two special cases.

$\bullet$ If the reduction following the expansion does not involve
identifying or omitting states of $\calA$, or equivalently if $\calA$
embeds in $\calB$, we say that $\calB$ is obtained from $\calA$ by a
\textit{reduced expansion} and we write $\calA \restep^{(p,w,q)} \calB$
or $\calA \restep \calB$.

$\bullet$ If the states $p$ and $q$ are equal to the distinguished
state $q_{0}$ of $\calA$, we say that $\calB$ is obtained from $\calA$
by an \textit{e-step} and we write $\calA \estep^{w} \calB$, or simply
$\calA \estep \calB$.

Finally, let $\calB$ be obtained from $\calA$ by identifying two
distinct vertices $p$ and $q$, and then reducing the resulting
automaton. Then we say that $\calB$ is obtained from $\calA$ by an
\textit{i-step} and we write $\calA \istep^{p = q} \calB$, or simply
$\calA \istep \calB$.

Note that if $\calA \expstep
\calB$, $\calA \restep \calB$, $\calA \estep \calB$ or $\calA \istep
\calB$, then $\calB$ is a reduced inverse automaton.

We first record a few facts.

\begin{fact}\label{fact1}
    Let $u$ be a reduced word labeling a path in $\calA$ from a state
    $p$ to a state $p'$, and from a state $q$ to a state $q'$,
    $$p \mapright{u} p', \quad q \mapright{u} q'.$$
    By definition of the reduction of dual automata, the identification
    of $p$ and $q$ implies that of $p'$ and $q'$, and the converse
    holds as well. Thus $\calA \istep^{p = q} \calB$ if and only if
    $\calA \istep^{p' = q'} \calB$.
\end{fact}    

Let us now examine in detail the effect of an operation of the form
$\expstep$.

\begin{fact}\label{fact2}
    Let $p, q$ be states of $\calA$ and let $w$ be a non-empty reduced
    word. Let $u$ be the longest prefix of $w$ that can be read in
    $\calA$ from state $p$, and let $v$ be the longest suffix of $w$
    that can be read in $\calA$ to state $q$ (that is, $v\inv$ is the
    longest prefix of $w\inv$ that can be read in $\calA$ from state
    $q$). We distinguish two cases:

    \begin{itemize}
	\item[(1)] If $|u| + |v| < |w|$, then $w = uw'v$ for some
	non-empty reduced word $w'$. If we let $p'$ (resp. $q'$) be the
	end (resp. start) state of the path labeled $u$ (resp. $v$) and
	starting in $p$ (resp. ending in $q$),
	$$p \mapright{u} p' \mapright{w'} q' \mapright{v} q,$$
	then the reduction process on the result of the expansion of
	$\calA$ by $(p,w,q)$ identifies the $|u|$ first edges and the
	$|v|$ last edges of the added path with existing edges of
	$\calA$, so that $\calA \expstep^{(p,w,q)} \calB$ if and only
	if $\calA \expstep^{(p',w',q')} \calB$ and the latter is a
	reduced expansion.
		
	\item[(2)] If $|u| + |v| \ge |w|$, then there exist words
	$x,y,z$, possibly empty, such that $u = xy$, $v = yz$ and $w =
	xyz$.  Let $p', p'', q', q''$ be the states of $\calA$ defined
	by the following paths
	$$p \mapright{x} p' \mapright{y} p'', \quad q' \mapright{y} q''
	\mapright{z} q.$$
	Then $\calA \expstep^{(p,w,q)} \calB$ if and only if $\calA
	\istep^{p' = q'} \calB$, if and only if $\calA \istep^{p'' =
	q''} \calB$.
    \end{itemize}    
\end{fact}

We derive from Fact~\ref{fact2} the following statement.

\begin{prpstn}\label{L04}
    Let $\calA$ and $\calB$ be inverse automata. If $\calA \estep^{w}
    \calB$, then $\calA \istep \calB$ or $\calA \restep^{(p,u,q)}
    \calB$ for some states $p$ and $q$ and a reduced word $u$ such
    that $|u| \le |w|$.
\end{prpstn}    

The following converse statements are derived from Facts~\ref{fact1}
and~\ref{fact2}.

\begin{prpstn}\label{H istep}
    Let $\calA$ be a reduced inverse automaton, let $H =
    L(\calA)\rho$, let $u$ and $v$ be reduced words labeling paths $q
    \mapright v q_{0} \mapright u p$ in $\calA$, and suppose that
    $\calA \istep^{p=q} \calB$.  Then $\calA \estep^{uv} \calB$ and
    $L(\calB)\rho = \langle H, uv\rangle$.
    
    In particular, $\rank(L(\calB)\rho) \le 1 + \rank(H)$.  If
    $\rank(L(\calB)\rho) = 1 + \rank(H)$, then $H \leff L(\calB)\rho$
    and if $C$ is a basis of $H$, then $C\cup\{uv\}$ is a basis of
    $L(\calB)\rho$.
\end{prpstn}

\begin{proof}
Let $\calA'$ be the expansion of $\calA$ by $(q_{0},uv,q_{0})$.  The
analysis in Fact~\ref{fact2} (2) shows that a step in the reduction of
$\calA'$ is provided by the automaton obtained in identifying $p$ and
$q$.  By Proposition~\ref{H(A)} we have $L(\calB)\rho = \langle H,
uv\rangle$ (hence the rank inequality), and the uniqueness statement
in Proposition~\ref{prop unique} then shows that $\calA \estep^{uv}
\calB$.

Let us now assume that $\rank(L(\calB)\rho) = 1 + \rank(H)$, and that
$C$ is a basis of $H$.  Then $C \cup \{uv\}$ is a generating set of
$L(\calB)\rho$ with cardinality equal to the rank of $L(\calB)\rho$,
so $C \cup \{uv\}$ is a basis of that subgroup by \cite[Prop.
I.3.5]{LS}.
\end{proof}

\begin{prpstn}\label{H restep}
    Let $\calA$ and $\calB$ be reduced inverse automata, let $w$ be a
    reduced word such that $\calA \restep^{(p,w,q)} \calB$, let $H =
    L(\calA)\rho$, and let $u$ and $v$ be reduced words labeling paths
    $q \mapright v q_{0} \mapright u p$ in $\calA$.  Then $\calA
    \estep^{uwv} \calB$ and $L(\calB)\rho = \langle H, uwv\rangle$.
    
    Moreover, $H \leff L(\calB)\rho$, $\rank(L(\calB)\rho) = 1 +
    \rank(H)$ and if $C$ is a basis of $H$, then $C \cup \{uwv\}$ is a
    basis of $L(\calB)\rho$.
\end{prpstn}

\begin{proof}
Since the expansion of $\calA$ by $(p,w,q)$ is a reduced expansion,
the word $uwv$ is reduced and the expansion by $(q_{0},uwv,q_{0})$
falls in the situation described in Fact~\ref{fact2}~(1).  Together
with Proposition~\ref{H(A)}, it follows that $\calA \estep^{uwv}
\calB$, which concludes the verification that of the first part of 
the proposition.

The free factor and the rank statements follow from the statement on 
a basis for $L(\calB)\rho$, which is a direct consequence of the 
definition of a reduced expansion and of the discussion on bases and 
spanning trees in Fact~\ref{end of sec 12}.
\end{proof}

We now introduce a measure of the \textit{cost} of a reduced
expansion or an i-step $\sigma$, written $\lambda(\sigma)$: if $\sigma$
is an i-step, then $\lambda(\sigma) = 0$; if $\sigma$ is a reduced
expansion, $\sigma = \restep^{(p,w,q)}$, its cost is the length of
$w$, $\lambda(\sigma) = |w|$. We extend this notion of cost to finite
sequences of i-steps and reduced expansions: if $\bar\sigma =
(\sigma_{1},\ldots,\sigma_{n})$ is such a sequence, we let
$$\lambda(\bar\sigma) =
(\lambda(\sigma_{1}),\ldots,\lambda(\sigma_{n})).$$
Finally, we introduce an order relation on the set of finite sequences
of non-negative integers. Let $\bar k = (k_{1},\ldots,k_{n})$ and
$\bar\ell = (\ell_{1},\ldots,\ell_{m})$ be such sequences. We say that
$\bar k \preceq \bar\ell$ if
\begin{itemize}
    \item[either] $n < m$,
    
    \item[or] $n = m$ and $\sum_{i=1}^n k_{i} < \sum_{i=1}^m
    \ell_{i}$,
    
    \item[or] $n = m$, $\sum_{i=1}^n k_{i} = \sum_{i=1}^m \ell_{i}$
    and $\bar k$ precedes $\bar\ell$ in the lexicographic order.
\end{itemize}    

It is routine to check that $\preceq$ is a well-order on the set of
finite sequences of non-negative integers, which is stable under the
concatenation of sequences. We write $\bar k\prec\bar\ell$ if $\bar
k\preceq\bar\ell$ and $\bar k\ne\bar\ell$.

\begin{prpstn}\label{L06}
    Let $\calA$, $\calA'$ and $\calB$ be inverse automata such that
    $\calA'$ is obtained from $\calA$ by a reduced expansion
    $\sigma_{1}$ and $\calB$ is obtained from $\calA'$ by an i-step
    $\sigma_{2}$,
    $$\calA \restep \calA' \istep \calB.$$
    Then there exist a sequence of reduced expansions or i-steps
    $\bar\sigma'$ of length 1 or 2 such that $\calB$ is obtained from
    $\calA$ by applying the steps in $\bar\sigma'$ and
    $\lambda(\bar\sigma') \prec \lambda(\sigma_{1},\sigma_{2})$.
\end{prpstn}

\begin{proof}
Suppose that $\calA \restep^{(p,w,q)} \calA' \istep^{r = s} \calB$
with $r\ne s$.  The cost of this sequence of transformations is
$(|w|,0)$.

Let $Q$ be the state set of $\calA$ and let $u$ and $v$ be reduced
paths,
$$q \mapright{v} q_0 \mapright{u} p.$$
Then $uwv$ is a reduced word and $L(\calA')\rho = \langle
L(\calA)\rho,uwv\rangle$ by Proposition~\ref{H restep}. We distinguish
three cases, depending whether or not $r$ and $s$ lie in $Q$.

\paragraph{Case 1: Both $r$ and $s$ are in $Q$.}
Let $x$ and $y$ be reduced words labeling paths in $\calA$ $s
\mapright{y} q_0 \mapright{x} r$.  Then the same words label similar
paths in $\calA'$ and it follows from Proposition~\ref{H istep} that
$$L(\calB)\rho = \langle L(\calA')\rho, xy\rangle = \langle
L(\calA)\rho, uwv, xy \rangle.$$

Let also $\calA''$ and $\calB'$ be determined by $\calA \istep^{r = s}
\calA'' \estep^{uwv} \calB'$. Then $L(\calB')\rho$ is also equal to
$\langle L(\calA)\rho, xy, uwv \rangle$, so that $\calB = \calB'$
by Proposition~\ref{prop unique}.

Note that the words $u$ and $v$ label paths from and into state $q_{0}$ in
$\calA''$ as well. It follows from Proposition~\ref{L04} that, if
$uwv \not\in L(\calA'')$, then $\calB$ can be obtained from
$\calA''$ by an i-step or by a reduced expansion of the form
$\restep^{(t,z,t')}$ with $|z|\le |w|$.

Thus $\calB$ is obtained from $\calA$ either by a sequence of 1 or 2
transformations, of cost $0$ or $(0,k)$ with $0\le k\le |w|$. This 
is $\prec$-less than $(|w|,0)$, as expected.

\paragraph{Case 2: Exactly one of $r$ and $s$ is in $Q$.}
Without loss of generality, we may assume that $r\in Q$ and $s\not\in
Q$.  Let $z$ be a reduced word labeling a path from $q_{0}$ to $r$ in
$\calA$, and hence also in $\calA'$.  Let $g$ be the unique reduced
word labeling a path from $p$ to $s$ in $\calA'$, using only edges
that were not in $\calA$.  By assumption, $g$ is a proper, non-empty
prefix of $w$.  Moreover, by Propositions~\ref{H(A)} and~\ref{H
istep},
$$L(\calB)\rho = \langle L(\calA')\rho, ugz\inv\rangle = \langle
L(\calA)\rho, uwv, ugz\inv \rangle.$$

Let $h$ be the longest common suffix of $g$ and $z$, so that $g = g'h$,
$z = z'h$, $g'{z'}\inv$ is reduced and we have the following paths in
$\calA'$,
$$q_0 \mapright{z'} r' \mapright{h} r, \quad p \mapright{g'} s'
\mapright{h} s.$$
Fact~\ref{fact1} shows that $\calA'\istep^{r' = s'} \calB$, so
we may assume that $h = 1$, $g = g'$ and $z = z'$. There is a
possibility that the word $g$ is now empty (if $h$ was in fact equal to
$g$), but in that case, we are returned to the situation of Case 1,
with $s' = p$. Thus we may still assume that $g\ne 1$. In particular,
the word $ugz\inv$ is reduced.

Then let $\calA''$ and $\calB'$ be defined by $\calA \estep^{ugz\inv}
\calA'' \estep^{uwv} \calB'$. Again $L(\calB')\rho = \langle
L(\calA)\rho, uwv, ugz\inv \rangle$, so $\calB = \calB'$ by
Proposition~\ref{prop unique}.

Proposition~\ref{L04} states that each e-step can be replaced by an
i-step or by a reduced expansion of cost bounded above by the cost of
the e-step.  Going back to Fact~\ref{fact2}, we see that the e-step
$\calA \estep^{ugz\inv} \calA''$ can be replaced by a transformation
of cost $k\le |g|$ since both $u$ and $z$ can be read from state
$q_{0}$ in $\calA$ (in fact, of cost exactly $|g|$ by definition of
$g$).  As for the e-step $\calA'' \estep^{uwv} \calB$, it can be
replaced by a transformation of cost $\ell \le |w|-|g|$ since $ug$ (a
prefix of $uw$) and $v$ can be read to state $q_{0}$ in $\calA''$.

Now, it suffices to verify that $(k,\ell) \prec (|w|,0)$, which is
easily done if we observe that $k + \ell \le |w|$ and $k < |w|$.

\paragraph{Case 3: Neither $r$ nor $s$ is in $Q$.}
Without loss of generality, we may assume that $r$ occurs before $s$
along the $w$-labeled path from $p$ to $q$. Thus, the word
$w$ factors as $w = w_{1}w_{2}w_{3}$ and the path in $\calA'$ made of
edges added to $\calA$ factors as
$$p \mapright{w_{1}} r \mapright{w_{2}} s \mapright{w_{3}} q.$$
Since $r\ne s$ and these vertices are not in $Q$, each of the three
factors $w_{1}, w_{2}, w_{3}$ is non-empty. Moreover,
$$L(\calB)\rho = \langle L(\calA')\rho, uw_{1}w_{3} v\rangle = \langle
L(\calA)\rho, uwv, uw_{1}w_{3} v \rangle.$$

Let $h$ be the longest common suffix of $w_{1}$ and $w_{3}\inv$, so
that $w_{1} = w'_{1}h$, $w_{3} = h\inv w'_{3}$, $w'_{1}w'_{3}$ is
reduced and we have the following paths in $\calA'$,
$$p \mapright{w'_{1}} r' \mapright{h} r \mapright{w_{2}} s
\mapleft{h} s' \mapright{w'_{3}} q.$$
Proposition~\ref{fact1} shows that $\calA'\istep^{r' = s'} \calB$, so
we may assume that $h = 1$, $w_{1} = w'_{1}$ and $w_{3} = w'_{3}$.
There is a possibility that the words $w_{1}$ or $w_{3}$ be now empty
(if $h$ was in fact equal to $w_{1}$ or $w_{3}$), but in that case, we
are returned to the situation of Cases 1 or 2, with $r' = p$ or $s' =
q$. Thus we may still assume that $w_{1}\ne 1$ and $w_{3}\ne 1$. In
particular, the word $uw_{1}w_{3}v$ is reduced.

Then let $\calA''$ and $\calB'$ be defined by $\calA
\estep^{uw_{1}w_{3}v} \calA'' \estep^{uwv} \calB'$. Then
$L(\calB')\rho = \langle L(\calA)\rho, uwv, uw_{1}w_{3}v \rangle$, so
$\calB = \calB'$ by Proposition~\ref{prop unique}.

As in Case 2, we use Fact~\ref{fact2} to verify that the e-step $\calA
\estep^{uw_{1}w_{3}v} \calA''$ can be replaced by a reduced expansion
of cost $k = |w_{1}w_{3}|$ since $u$ and $v$ are the maximal prefix
and suffix of $uw_{1}w_{3}v$ that can be read from and to state
$q_{0}$ in $\calA$.  As for the e-step $\calA'' \estep^{uwv} \calB$,
it can be replaced by a reduced expansion of cost $\ell = |w_{2}|$
since $uw_{1}$ and $w_{3}v$ are the maximal prefix and suffix of $uwv$
that can be read from and to state $q_{0}$ in $\calA''$.

Now, it suffices to verify that $(k,\ell) \prec (|w|,0)$, which is
easily done if we observe that $k + \ell = |w|$ and $k < |w|$.
\end{proof}

%%%%%%%%%%%%%%%%%%%%%%%%%%%%%%%
\section{Deciding the free factor relation}\label{sec deciding}

%%%%%%%%%%%%%%%%%%%%%%%%%%%%%%%
\subsection{A geometric characterization of free factors}

We put together the technical results from Section~\ref{careful} to
prove the following characterization of free factors.

\begin{thrm}\label{charact ff}
    Let $H, K$ be finitely generated subgroups of $F = F(A)$ and assume
    that $d = \rank(K) - \rank(H) > 0$. Then $H$ is a free factor of
    $K$ if and only if the inverse automaton $\Gamma_{A}(H)$ can be
    transformed in $\Gamma_{A}(K)$ by a sequence of $d'\le d$ i-steps
    followed by $d - d'$ reduced expansions.
\end{thrm}

\begin{proof}
We first observe that $H$ is a free factor of $K$ if and only if there
exist $d$ elements $k_{1},\ldots,k_{d}$ of $F(A)$ such that $\langle
H\cup\{k_{1},\ldots,k_{d}\}\rangle = K$. This follows from the fact
that an $r$-element generating set in a rank $r$ free group, is a basis
\cite[Prop. I.3.5]{LS}.

By definition of e-steps, this means that $H\leff K$ if and only if
$\Gamma_{A}(H)$ yields $\Gamma_{A}(K)$ by a sequence of $d$ e-steps.

Now Propositions~\ref{L04}, \ref{H istep} and~\ref{H restep} show that
this is equivalent to the fact that $\Gamma_{A}(H)$ yields
$\Gamma_{A}(K)$ by a sequence of $d$ i-steps or reduced expansions.

Since $\preceq$ is a well-order on the set of finite sequences of
non-negative integers, we may consider a sequence $\bar\sigma$ of $d$
i-steps and reduced expansions leading from $\Gamma_{A}(H)$ to
$\Gamma_{A}(K)$, which is $\preceq$-minimal.  Proposition~\ref{L06},
together with the stability of $\preceq$ under concatenation, then
shows that the i-steps in $\bar\sigma$ come before the reduced
expansions.  Thus, $H\leff K$ if and only if $\Gamma_{A}(H)$ yields
$\Gamma_{A}(K)$ by a sequence of $d'$ i-steps followed by $d - d'$
reduced expansions.
\end{proof}

It follows from the discussion on bases and spanning trees in
Fact~\ref{end of sec 12} that if $\Gamma_{A}(H)$ embeds in
$\Gamma_{A}(K)$, then $H$ is a free factor of $K$.  Not every free
factor of $K$ occurs that way, and those that do are called the
\textit{graphical free factors of $K$} (\textit{with respect to $A$}).
It is easily verified that $H$ is a graphical free factor of $K$ if
and only if $\Gamma_{A}(H)$ yields $\Gamma_{A}(K)$ by a sequence of
reduced expansions.

\begin{crllr}\label{coroll HffK}
    Let $H, K$ be finitely generated subgroups of $F = F(A)$ and assume
    that $d = \rank(K) - \rank(H) > 0$. Then $H$ is a free factor of
    $K$ if and only if the inverse automaton $\Gamma_{A}(H)$ can be
    transformed by a sequence of $d'\le d$ i-steps into some 
    $\Gamma_{A}(L)$ such that $\rank(L) = \rank(H) + d'$ and $L$ is a 
    graphical free factor of $K$ with respect to $A$.
\end{crllr}

\begin{proof}
Let us first assume that $H\leff K$.  By Theorem~\ref{charact ff}, for
some $d'\le d$, $\Gamma_{A}(H)$ can be taken to some $\Gamma_{A}(L)$
by a sequence of $d'$ i-steps, and $\Gamma_{A}(L)$ can be taken to
$\Gamma_{A}(K)$ by a sequence of $d - d'$ re-steps.  Since $\rank(K) =
\rank(H) + d$ and an i-step or an re-step can increment the rank by at
most one (Propositions~\ref{H istep} and~\ref{H restep}), $L$ must have
rank exactly $\rank(H) + d'$.

Conversely, suppose that a sequence of $d' \le d$ i-steps takes
$\Gamma_{A}(H)$ to $\Gamma_{A}(L)$ in such a way that $\rank(L) =
\rank(H) + d'$ and $\Gamma_{A}(L)$ embeds in $\Gamma_{A}(K)$.  By
Propositions~\ref{H istep} and~\ref{H restep} again, we have $H \leff
L \leff K$, and hence $H\leff K$.
\end{proof}

\begin{rmrk}\label{rank incrementing}
    We observe the following by-product of the proof of
    Corollary~\ref{coroll HffK}: if $\Gamma_{A}(H)$ can be transformed
    into $\Gamma_{A}(L)$ by a sequence of $d'$ i-steps such that
    $\rank(L) = d' + \rank(H)$, then for every i-step $\calA \istep
    \calB$ occurring in that sequence, we have $\rank(L(\calB)\rho) =
    1 + \rank(L(\calA)\rho)$.  We say that such an i-step is
    \textit{rank-incrementing}.
\end{rmrk}

In the special case where $K = F$, we have the following statement.

\begin{crllr}\label{coroll HffF}
    Let $H$ be a finitely generated subgroup of $F = F(A)$, let
    $A_{0}$ be the set of letters in $A$ that occur in $\Gamma_{A}(H)$
    and let $d = |A_{0}| - \rank(H) = \rank(F(A_{0})) - \rank(H)$.
    Then $H$ is a free factor of $F$ if and only if $d \ge 0$ and
    $\Gamma_{A}(H)$ can be transformed into a one-vertex automaton by
    a sequence of $d$ i-steps.
\end{crllr}

\begin{proof}
By Corollary~\ref{coroll HffK}, $H$ is a free factor of $F(A)$ if and
only if $\Gamma_{A}(H)$ yields a one-vertex automaton $\calB$ by a
sequence of $d' \le |A| - \rank(H)$ i-steps, in such a way that
$L(\calB)\rho$ has rank $d' + \rank(H)$.  Now the set of letters
occurring in such an automaton $\calB$ must be exactly $A_{0}$, so
$L(\calB)\rho = F(A_{0})$ and $d' = d$.
\end{proof}

%%%%%%%%%%%%%%%%%%%%%%%%%%%%%%%
\subsection{Deciding whether $H\leff F$}\label{sec algorithm}

We first consider the problem of deciding whether a given finitely
generated subgroup $H$ of $F = F(A)$ is a free factor of $F$.  With
the notation of Corollary~\ref{coroll HffF} and Remark~\ref{rank
incrementing}, the algorithm to decide whether $H\leff F$ consists of
the following.  We need to explore the sequences of rank-incrementing
i-steps, starting from $\Gamma_{A}(H)$ and of length $d = |A_{0}| -
\rank(H)$.  Then $H\leff F$ if and only if one of the automata
occurring at the end of one of these sequences has a single vertex. 
Note that no automaton obtained after less than $d$ i-steps could be 
a single-vertex automaton.

These automata can be viewed as nodes of a tree, rooted at
$\Gamma_{A}(H)$, in which the children of a node are the automata
produced by a rank-incrementing i-step.  Let $n$ be the number of
states of $\Gamma_{A}(H)$.  Then every automaton occurring along this
tree has at most $n$ states.

If $\calA$ is such an automaton, then $\calA$ has at most
$\frac12(n^2-n)$ pairs of distinct vertices, and hence at most
$\frac12(n^2-n)$ children, each of which has at most $n - 1$ states.
The computation of the children of $\calA$ is done by computing all
the (at most $\frac12(n^2-n)$) possible i-steps, computing
the ranks of the corresponding subgroups and retaining only those of
rank $1 + \rank(L(\calA)\rho)$.  It follows from Fact~\ref{end of sec
12} that the cost of the computation of the children of $\calA$ is
$O(n^{4})$.

Moreover, there are at most $O(n^{2d-2})$ nodes of the tree at depth
at most $d-1$, and the computation of these nodes and their children
requires time $O(n^{2d+2})$.  For each of the $O(n^{2d})$ automata at
depth $d$, the verification whether the automaton has a single node
takes constant time, so the total cost of the algorithm is
$O(n^{2d+2})$.

Finally, if $H$ is given by a finite set of generators, of total
length $\ell$, we recall that computing $\Gamma_A(H)$ takes time
$O(\ell^2)$ and that $\Gamma_A(H)$ has at most $\ell$ states and
$\ell$ edges.  This discussion justifies the following statement.

\begin{thrm}\label{thm35}
    Let $A$ be a fixed alphabet.  Then there is an algorithm which,
    given reduced words $h_{1},\ldots,h_{n}$ in $F(A)$ with total
    length $\ell$, decides whether the subgroup $H$ generated by the
    $h_{i}$ is a free factor of $F(A)$ in time $O(\ell^{2d+2})$, where
    $d = |A_{0}| - \rank(H)$ and $A_{0}$ is the set of letters in $A$
    that occur in the $h_{i}$.
\end{thrm}    

\begin{rmrk}
    The tree exploration described above can be speeded up by the
    following observation.  If the automaton $\calA$ occurs in a
    sequence of rank-incre\-menting i-steps from $\Gamma_{A}(H)$ to the
    one-vertex automaton $\Gamma_{A}(F(A_{0}))$ (a \textit{winning
    sequence}), then $L(\calA)\rho \leff F(A_{0})$, which implies that
    $L(\calA)\rho$ is a free factor of every subgroup of $F(A_{0})$
    containing it.  It follows that, if some i-step $\calA \istep
    \calB$ does not increment the rank, then $\calA$ does not occur in
    a winning sequence, that is, we may ignore the subtree below
    $\calA$.
    
    There are undoubtedly other implementation tricks and ideas that
    can reduce the decision process, however without changing the
    worst-case complexity.
\end{rmrk}

\begin{rmrk}
    In view of Touikan's announcement (see Fact~\ref{end of sec 12}),
    the time complexity in the above theorem can be lowered to
    $O(\ell^{2d+1}\log^*\ell)$.
\end{rmrk}

The above discussion of complexity depends on a \textit{uniform cost
assumption}, by which we assume that the elementary operations on $A$
(reading or writing a letter, comparing two letters) require unit
time.  In other words, we have been assuming that $A$ is fixed, and
not part of the input.  We will shortly consider the problem of
deciding whether $H \leff K$ where both $H$ and $K$ are part of the
input, and for the purpose of that discussion, we reconsider the
algorithm and the complexity established above under the \textit{bit
cost assumption}: we let $r$ be the cardinality of $A$, and we
consider that each letter is identified by a bit string of length at
most $\log r$, so that the elementary operations on $A$ require
$O(\log r)$ units of time.

Let $f(\ell,r)$ be the complexity of reducing a dual automaton on
$\tilde A$ with at most $\ell$ vertices and $\ell$ edges.  The
computation of the automaton obtained by an i-step from an
$\ell$-vertex automaton takes time at most $f(\ell,r)$.  To compute
the tree of rank-incrementing i-steps rooted at $\Gamma_A(H)$, we need
to compute the at most $O(\ell^2)$ children of at most $\ell^{2d-2}$
nodes, which requires time $O(\ell^{2d}f(\ell,r))$.  Finally, deciding
whether an automaton has a single vertex takes constant time, so the
total complexity of the algorithm is $O(\ell^{2d}f(\ell,r))$.

We now give an upper bound of $f(\ell,r)$.  Let $\calA$ be an
$A$-automaton with at most $\ell$ states and $\ell$ edges.  The
identifiers of states require space $O(\log\ell)$ and the identifiers
of letters require space $O(\log r)$.  Moreover, we assume that the
set of states and the set $\tilde A$ are linearly ordered, and
equipped with a constant time \textsf{next} function.

The automaton can be viewed as the lexicographically ordered list of
all triples $(u,a,v)$ such that either $a\in A$ and there is an
$a$-labeled edge from state $u$ to state $v$, or $\bar a\in A$ and
there is a $\bar a$-labeled edge from state $v$ to state $u$.  Each
entry of this list requires space $O(\log\ell + \log r) = O(\log(\ell
r))$, and the list contains at most $\ell$ entries.  In particular, a
complete scan of the list takes time $O(\ell \log(\ell r))$.

The reduction of $\calA$ consists in performing elementary reductions
of type 1 as long as it is possible, and then elementary reductions of
type 2.  To find out whether a type 1 reduction is possible, one needs
to scan the list to find two (consecutive) cells of the form $(u,a,v)$
and $(u,a,v')$, requiring $O(\ell \log(\ell r))$.  Performing the
identification consists in reading through the list, replacing every
occurrence of $v'$ by $v$, and reordering the list: this takes time
$O(\ell \log(\ell r))$.  To find out whether an elementary reduction
of type 2 is possible, one scans the list to find a vertex $u$ such
that there is a unique triple of the form $(u,a,v)$ in the list:
again, this takes time $O(\ell \log(\ell r))$.  Finally, performing
the reduction consists simply in deleting the entries $(u,a,v)$ and
$(v,\bar a,u)$ in the list.

Summarizing, since we will perform at most $\ell$ reductions, we can
take the function $f(\ell,r)$ to be equal to $\ell^2\log(\ell r)$. 
This yields the following statement.

\begin{thrm}\label{thm-generalized35}
    There is an algorithm which, given an alphabet $A$ of size $r$ and
    reduced words $h_{1},\ldots,h_{n} \in F(A)$ with total length
    $\ell$, decides whether the subgroup $H$ generated by the $h_{i}$
    is a free factor of $F(A)$ in time $O(\ell^{2d+2} \log(\ell r))$,
    where $d = |A_{0}| - \rank(H)$ and $A_{0}$ is the set of letters
    in $A$ that occur in the $h_{i}$.
\end{thrm}

%%%%%%%%%%%%%%%%%%%%%%%%%%%%%%%
\subsection{Deciding whether $H \leff K$}\label{sec algorithm HK}

We now suppose that $H$ and $K$ are subgroups of $F(A)$, given by sets
of generators with total length $\ell$, and we consider the problem of
deciding whether $H\leff K$.  The algorithm is the following: we first
compute $\Gamma_{A}(K)$ and we verify whether each generator of $H$ 
can be read as a loop at $q_0$, the designated vertex of 
$\Gamma_A(K)$. If not, then $H$ is not contained in $K$, and cannot 
be a free factor.

We now assume that $H \le K$ and we compute a spanning tree $T$ of the
graph $\Gamma_{A}(K)$.  This spanning tree determines a basis $B$ of
$K$, indexed by the edges of $\Gamma_{A}(K)$ that are not in $T$ (see
Fact~\ref{end of sec 12}), and we express the given generators of $H$
in terms of this basis: it suffices to read the generators of $H$ as
loops at the designated vertex in $\Gamma_{A}(K)$, and to record the
sequence of edges traversed and not in $T$.  In particular, each
generator of $H$ is expressed as a reduced word in $\tilde B^*$, that
is no longer than its expression as a reduced word in $\tilde A^{*}$.
We then use this expression of the generators of $H$ to construct
$\Gamma_{B}(H)$, and the algorithm in Section~\ref{sec algorithm} to
decide whether $H$ is a free factor of the ambient free group, namely
$K$.

Let us now discuss the complexity of this algorithm, assuming that 
$H$ and $K$ are given by tuples of generators of total length $\ell$ 
and that $F$ has rank $r$. As discussed in the previous section, 
computing $\Gamma_A(K)$ takes time $O(\ell^2\log(\ell r))$. Running 
a word $w \in F(A)$ in the automaton $\Gamma_A(K)$ requires reading 
sequentially each letter of $w$ (time $O(\log r)$) and looking for 
the corresponding transition in the table representing the automaton 
(time $O(\ell\log(\ell r))$). If $w$ has length  $\ell'$, this takes time 
$O(\ell'\ell\log(\ell r))$. In order to decide whether $H$ is contained 
in $K$, this must be done for every generator of $H$, and since the 
total length of these generators is at most $\ell$, this requires 
time $O(\ell^2\log(\ell r))$.

Assuming that $H$ is indeed contained in $K$, the next step is to
construct a spanning tree $T$ of $\Gamma_A(H)$, for instance by
marking certain edges in the list representing $\Gamma_A(K)$.  Again,
this can be done in time $O(\ell^2\log(\ell r))$.  The positively
labeled edges of $\Gamma_A(K)$ not in $T$ are in bijection with a
basis $B$ of $K$.  In particular, the rank of $K$ is at most $\ell$,
an upper bound of the number of edges in $\Gamma_A(K)$.  The elements
of $B$, seen as words in $F(A)$, consist of the label of a path in
$T$, followed by the label of an edge not in $T$, followed by the
label of a path in $T$.  In particular, their length is at most twice
the number of vertices of $\Gamma_A(K)$ plus one, that is $O(\ell)$.
But we do not need to compute these words: it suffices to number (from
$1$ to $\rank(K) \le \ell$) the positively labeled edges not in $T$.
Then, for each generator $h$ of $H$, reading $h$ in $\Gamma_A(K)$ from
state $q_0$ and keeping track of the (identifier of the) edges
traversed that are not in $T$, provides an expression of $h$ in $B$.
Moreover, the total length of the expression of the generators of
$H$ in this basis of $K$ is at most $\ell$.

We now simply apply Theorem~\ref{thm-generalized35} to a length $\ell$
set of generators of a subgroup of a free group of rank at most
$\ell$, to justify the following corollary.

\begin{crllr}
    Given tuples of generators for subgroups $H$ and $K$ of $F(A)$,
    with total length $\ell$, one can decide whether the subgroup $H$
    is a free factor of $K$ in time
    $O(\ell^{2d+2}\log(\ell r))$, where $d = \rank(K) - \rank(H)$.
\end{crllr}

\begin{rmrk}\label{alternative algo H leff K}
    Corollary~\ref{coroll HffK} suggests an alternative algorithm to
    decide whether $H\leff K$: one may explore the sequences of
    rank-incrementing i-steps of length at most $d = \rank(K) -
    \rank(H)$, starting from $\Gamma_{A}(H)$ and producing only
    representations of subgroups contained in $K$.  Each
    $\Gamma_{A}(L)$ occurring in such a sequence must be checked to
    verify whether it can be embedded in $\Gamma_{A}(K)$: if it can,
    then $H\leff K$, and if it cannot, then the automata produced by a
    rank-incrementing i-step from it must be computed and checked in
    their turn.
%     
%     This alternative algorithm is more complex to implement, and it 
%     does not help in terms of complexity.
\end{rmrk}

%%%%%%%%%%%%%%%%%%%%%%%%%%%%%%%
\subsection{Complement of a free factor}

By definition of a free factor (see the introduction), $H\leff K$ if
there exists a basis $C$ of $H$ and a disjoint set $D\subseteq K$ such
that $C \cup D$ is a basis of $K$.  In that case, the subgroup
generated by $D$ is called a \textit{complement of $H$ in $K$}.  It is
easy to see that this complement is not uniquely defined.

If $H$ is indeed a free factor of the free group $F$, the algorithm
described in Section~\ref{sec algorithm} also provides a sequence of
rank-incrementing i-steps taking $\Gamma_{A}(H)$ to
$\Gamma_{A}(F(A_{0}))$, where $A_{0}$ is the set of letters of $A$
that occur in the generators of $H$.  Repeated applications of
Proposition~\ref{H istep} then yield a basis of $F(A_{0})$ containing
a basis of $H$, and hence a complement of $H$ in $F(A_{0})$.

We can be a little more precise.  If $\calA$ is a reduced inverse
automaton, let $\diam_{q_{0}}(\calA)$, the \textit{$q_{0}$-diameter}
of $\calA$, be the longest shortest path in $\calA$ from $q_{0}$ to
a state, that is
$$\diam_{q_{0}}(\calA) = \max_{q}\Big(\min\{|u| \mid \textrm{$u$ is a
reduced word and $q_{0} \mapright u q$ in $\calA$}\}\Big).$$
Then Proposition~\ref{H istep} shows that if $\calA \istep \calB$ and 
$L(\calA)\rho \leff L(\calB)\rho$, then a complement of 
$L(\calA)\rho$ is generated by a word of length at most 
$2\ \diam_{q_{0}}(\calA)$.

This justifies the following statement.

\begin{prpstn}
    Let $H$ be a free factor of $F(A)$.  Then one can construct
    effectively a basis for a complement of $H$ in $F(A)$, consisting
    of words of length at most $2\ \diam_{q_{0}} (\Gamma_{A}(H))$.
\end{prpstn}

It is interesting to observe that this is a graphical analogue -- and
a minor improvement -- of Federer and J\'onsson's result (see
\cite[Prop.  I.2.26]{LS}) mentioned in Section~\ref{whitehead}.
Indeed, if $h_{1},\ldots, h_{n}$ are a set of generators of $H$, then
every edge of $\Gamma_{A}(H)$ is on a loop at $q_{0}$ labeled by some
$h_{i}$, and hence $\max_{i}|h_{i}| \ge 2\
\diam_{q_{0}}(\Gamma_{A}(H))$.

Let us now consider the problem of constructing a basis of a
complement of $H$ in $K$, where $H$ and $K$ are given finitely
generated subgroups of $F(A)$ and $H\leff K$.  The algorithm described
in Remark~\ref{alternative algo H leff K} (based on
Corollary~\ref{coroll HffK}) shows that one can construct effectively
a graphical free factor $L$ of $K$ such that $H \leff L$.  As above,
in view of Proposition~\ref{H istep}, the same algorithm can produce a
basis of a complement of $H$ in $L$ consisting of words of length at
most $2\ \diam_{q_{0}}(\Gamma_{A}(H))$.  There remains to construct a
basis for a complement of $L$ in $K$.

In view of Fact~\ref{end of sec 12} (and using the notation therein),
one can find a (basis of a) complement of $L$ in $K$ by considering a
spanning tree of $\Gamma_{A}(L)$, extending it to a spanning tree $T$
of $\Gamma_{A}(K)$, and considering the elements $b_{e}$ associated
with edges of $\Gamma_{A}(K)$ that are neither in $\Gamma_{A}(L)$ nor
in $T$. In particular, the words $b_{e}$ have length at most $1 + 2\ 
\diam_{q_{0}}(\Gamma(K))$. Thus we have the following statement.

\begin{prpstn}
    Let $H, K$ be finitely generated subgroups of $F(A)$.  If $H\leff
    K$, then one can construct effectively a basis for a complement of
    $H$ in $K$, consisting of words of length at most $\max(2\
    \diam_{q_{0}}(\Gamma_{A}(H)), 1 + 2\
    \diam_{q_{0}}(\Gamma_{A}(K)))$.
\end{prpstn}

%%%%%%%%%%%%%%%%%%%
\begin{acknowledgement}
The authors wish to thank both anonymous referees for their many
helpful comments, which have substantially contributed to the
improvement of the paper.  All remaining shortcomings are naturally
ours.
\end{acknowledgement}

%%%%%%%%%%%%%%%%%%%%%%%%

\end{document}